\documentclass[11pt,twoside,reqno]{amsart}
\usepackage{graphicx}
\usepackage{amsmath}
\usepackage{amsfonts}
\usepackage{amssymb}
\usepackage{amscd}
\usepackage{amsmath}
\usepackage{amscd}
\usepackage{amsthm}
\usepackage{amssymb}
\usepackage{amsfonts}
\usepackage{epsfig}
\usepackage{amsmath,amscd}
\usepackage{graphicx}
\usepackage{geometry,enumerate,latexsym,tabularx,shapepar}
\usepackage{pslatex}
\usepackage[all,2cell,knot,poly]{xy} \UseAllTwocells \SilentMatrices
\usepackage{slashbox}
\input{xypic}
\xyoption{all}
\numberwithin{equation}{section}
\newcommand{\ie}{{\it i.e.\/}\ }

\newtheorem{thm}{Theorem}[section]

\newtheorem{lem}[thm]{Lemma}
\newtheorem{prop}[thm]{Proposition}
\theoremstyle{definition}
\newtheorem{defn}[thm]{Definition}
\theoremstyle{remark}
\newtheorem{rem}[thm]{Remark}

\numberwithin{equation}{section}
\newtheorem{exa}[thm]{Example}

\newcommand{\bq}{\begin{eqnarray}}
\newcommand{\nq}{\end{eqnarray}}

\renewcommand{\d}{\partial}

\def\sign{{\rm sign}}

\newcommand{\R}{\mathbb{R}}
\newcommand{\s}{\mathbf{S^3}}

\def\R{{\mathbb R}}

\def\cV{{\mathcal V}}

\def\sign{{\rm sign}}

\begin{document}

\title[Khovanov Homology and embedded Graphs ]{Khovanov Homology and embedded Graphs }%
\author{Ahmad Zainy Al-Yasry} %

\address{University of Baghdad\\ and\\Abdus Salam International Center for Theoretical Physics (ICTP)}%
\email{ahmad@azainy.com}%


\date{2009}%


\begin{abstract}
We construct a cobordism group for embedded graphs in two different ways, first by using  sequences
of two basic operations, called ``fusion" and ``fission", which in terms of
cobordisms correspond to the basic cobordisms obtained by attaching or
removing a 1-handle, and the other one by using the concept of a 2-complex
surface with boundary is the union of these knots. A discussion given to the question
 of extending Khovanov homology from links to embedded graphs, by using the
 Kauffman topological invariant of embedded graphs by associating family of links and
 knots  to a such graph by using some local replacements at each vertex in the graph.
 This new concept of Khovanov homology of an embedded graph constructed to be the sum of the
Khovanov homologies of all the links and knots.
\end{abstract}

\maketitle

\section{Introduction}

The construction of the cobordism group for links and for knots and their
relation is given in \cite{Hoso}. We then consider the question of constructing a
similar cobordism group for embedded graphs in the 3-sphere.
We show that this can actually be done in two different ways, both
of which reduce to the same notion for links. The first one
comes from the description of the cobordisms for links in terms of sequences
of two basic operations, called ``fusion" and ``fission", which in terms of
cobordisms correspond to the basic cobordisms obtained by attaching or
removing a 1-handle. We define analogous operations of fusion and fission
for embedded graphs and we introduce an equivalence relation of cobordism
by iterated application of these two operations.\\ The second possible definition
of cobordism of embedded graphs is a surface (meaning here 2-complex)
in $S^3\times [0,1]$ with boundary the union of the given graphs.
 While for links, where cobordisms are realized
by smooth surfaces, these can always be decomposed into a sequence of handle
attachments, hence into a sequence of fusions and fissions, in the case of
graphs not all cobordisms realized by 2-complexes can be decomposed
as fusions and fissions, hence the two notions are no longer equivalent.
The idea of categorification the Jones  polynomial is known by
Khovanov Homology for links which is a new link invariant introduced by
Khovanov \cite{kh},\cite{Ba}. For each link $L$ in $S^3$ he defined a graded chain complex,
with grading preserving differentials, whose graded Euler characteristic is
equal to the Jones polynomial of the link $L$. The idea of Khovanov Homology
for graphs arises from the same idea of Khovanov homology for links by
the categorifications the chromatic polynomial of graphs.
This was done by L. Helme-Guizon and Y. Rong \cite{laur},
for each graph G, they defined a graded chain complex whose graded
Euler characteristic is equal to the chromatic polynomial of G.
In our work we try to recall, the Khovanov homology for links and graphs.\\
We discuss the question of extending Khovanov homology from links to
embedded graphs. This is based on a result of Kauffman that constructs a topological
invariant of embedded graphs in the 3-sphere by associating to such a graph
a family of links and knots obtained using some local replacements at
each vertex in the graph. He showed that it is a topological invariant by showing
that the resulting knot and link types in the family thus constructed are invariant
under a set of Reidemeister moves for embedded graphs that determine the
ambient isotopy class of the embedded graphs. We build on this idea and simply
define the Khovanov homology of an embedded graph to be the sum of the
Khovanov homologies of all the links and knots in the Kauffman invariant
associated to this graph. Since this family of links and knots is a topologically invariant,
so is the Khovanov homology of embedded graphs defined in this manner. We close this paper
by giving an example of computation of Khovanov homology for an embedded graph
using this definition.\\
\\
\textbf{\textrm{Acknowledgements:}} The author would like to express his deeply grateful to Prof.Matilde Marcolli
for her advices and numerous fruitful discussions. He is also thankful to Louis H. Kauffman and Mikhail Khovanov
for their support words and advices. Some parts of this paper done in Max-Planck Institut f\"ur Mathematik (MPIM), Bonn, Germany,
the author is kindly would like to thank MPIM for their hosting and subsidy during his study there.


\section{Knots and Links Cobordism Groups}
A notion of knot cobordism group and link cobordism group can be given by using cobordism classes of knots and links
to form a group \cite{foxmi},\cite{Hoso}. The link cobordism group splits into the direct sum of the knot cobordism group
and an infinite cyclic group which  represents the linking number, which is invariant under link cobordism \cite{Hoso}.
In this part we will give a survey about both  knot and link cobordism groups. In a later part of
this work we will show that the same idea can be adapted to construct
a graph cobordism group.
\\
\subsection{Knot cobordism group}
We recall the concept of cobordism between knots introduced in \cite{foxmi}. Two knots $K_1$ and $K_2$ are called knot cobordic
if there is a locally flat cylinder $ S $ in $\s \times [0,1]$  with boundary $\partial S =K_1 \cup -K_2$ where
$ K_1 \subset \s \times \{0\} $ and $ K_2 \subset \s \times \{1\} $. We then  write $ K_1 \sim K_2$.
The critical points in the cylinder are assumed be  minima (birth), maxima (death), and saddle points.
In the birth point at some $t_0$ there is a sudden appearance of a point. The point becomes an unknotted circle in the level immediately above $t_0$.
At a maxima or death point, a circle collapses to a point and disappearance from higher levels.\\
\\
\begin{figure}[htp]
\begin{center}
\includegraphics[width=10cm]{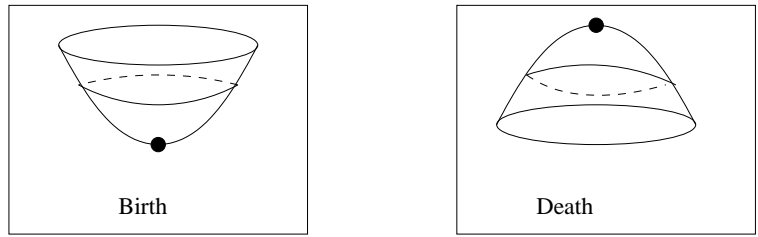}
\caption{ Death and Birth Points }\label{fig12}
\end{center}
\end{figure}
\\
For the saddle point, two curves touch and rejoin as illustrated in figure (\ref{fig13})
\\
\begin{figure}[htp]
\begin{center}
\includegraphics[width=10cm]{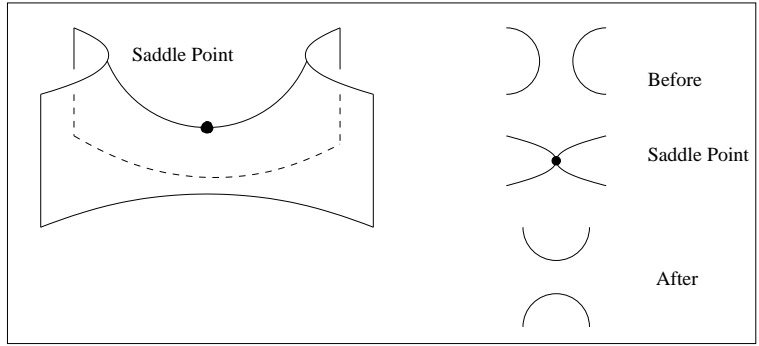}
\caption{ Saddle Point}\label{fig13}
\end{center}
\end{figure}
\\
  These saddle points are of two types:
{\em negative} if with increasing $t$  the number of the cross sections decreases and {\em positive}
if the number increases.\\ A transformation from one cross section to another is called negative
hyperbolic  transformation if there is only one saddle point between the two
cross sections and if the number of components decreases. We can define analogously a positive hyperbolic
transformation.
\begin{defn}\cite{Hoso}
We say that two knots $K_1$ and $K_2$ are related by an elementary cobordism
if the knot $K_2$ is obtained by $r-1$ negative hyperbolic transformations from a split link consisting of $K_1$
together with $r-1$ circles.
\end{defn}
What we mean by split link is a link with $n$ components $(K_i, i=1....n)$ in $\s$ such that there are mutually disjoint
$n$ 3-cells $(D_i, i=1....n)$ containing $K_i, i=1,2...,n $
\begin{lem}\cite{Hoso}
Two knots are called knot cobordic if and only if they are related by a sequence of elementary cobordisms
\end{lem}
It is well known that the oriented knots form a commutative semigroup under the operation of composition $\#$.
Given two knots $ K_1$ and $K_2$, we can obtain a new knot by removing a small arc from each knot and then connecting
the four endpoints by two new arcs. The resulting knot is called the composition of the two knots
$ K_1$ and $K_2$  and is denoted by $ K_1 \# K_2$.\\
Notice that, if we take the composition of a knot $K$ with the unknot $\bigcirc$ then the result is again $K$.
\begin{lem}
The set of oriented knots with the connecting operation $\#$ forms a semigroup with identity  $\bigcirc$
\end{lem}
Fox and Milnor \cite{foxmi} showed that composition of knots induces a composition on knot cobordism classes $[K] \# [K']$.
This gives an abelian group $G_K$ with $[\bigcirc]$ as identity and the negative is $ -[K]= [-K]$,
where the $-K$ denotes the reflected inverse of $K$.
\begin{thm}
The knot cobordism classes with the connected sum operation $\#$ form an abelian group, called the knot cobordism
group and denoted by $ G_K$.
\end{thm}

\subsection{Link cobordism group}\cite{Hoso}
For links, the conjunction operation $\&$ between two links gives a commutative semigroup. $ L_1 \& L_2 $
is a link represented by the union of the two links $ l_1 \cup l_2 $ where $ l_1$ represents $ L_1$ and $ l_2$
represents $ L_2$ with mutually disjoint 3-cells $ D_1$ and $ D_2$ contain $ l_1$ and $ l_2$ respectively.
Here ``represents" means that we are working with ambient isotopy classes $L_i$ of links (also
called {\em link types}) and the $l_i$ are chosen representatives of these classes. In the following
we loosely refer to the classes $L_i$ also simply as links, with the ambient isotopy equivalence
implicitly understood.
The zero
of this semigroup is the link consisting of just the empty link. The link cobordism group is constructed using the conjunction operation
and the cobordism classes. We recall below the definition of cobordism of links.\\ Let $L$ be a link in $\s$  containing $r$-components
$ k_1,....,k_r$ with a split sublink $ L_s =k_1 \cup k_2 \cup ....\cup k_t $, $ t \leq r$ of $L$.
Define a knot $ {\bf{\hat{K}}}$ to be $ k_1 + k_2 + ....+ k_t + \partial B_{t+1}+ \partial B_{t+2}+....+ \partial B_{r}$ where
$\{ B_{t+1}, B_{t+2}, B_{t+3}....,B_{r}\}$  are disjoint bands in $\s$  spanning $L_s$  \cite{Hoso}. The operation $+$
means additions in the homology sense. Put $ L_1 = L_s \cup  k_{t+1} \cup k_{t+2}....\cup k_{r}$ and
$L_2 = {\bf{\hat{K}}} \cup k_{t+1} \cup k_{t+2}....\cup k_{r}$. Now, the operation of replacing $L_1$ by $L_2$ is called
{\em{\bf{fusion}}} and $L_2$ by $ L_1$ is called {\em{\bf{fission}}}.\\
\begin{figure}[htp]
\begin{center}
\includegraphics[width=10cm]{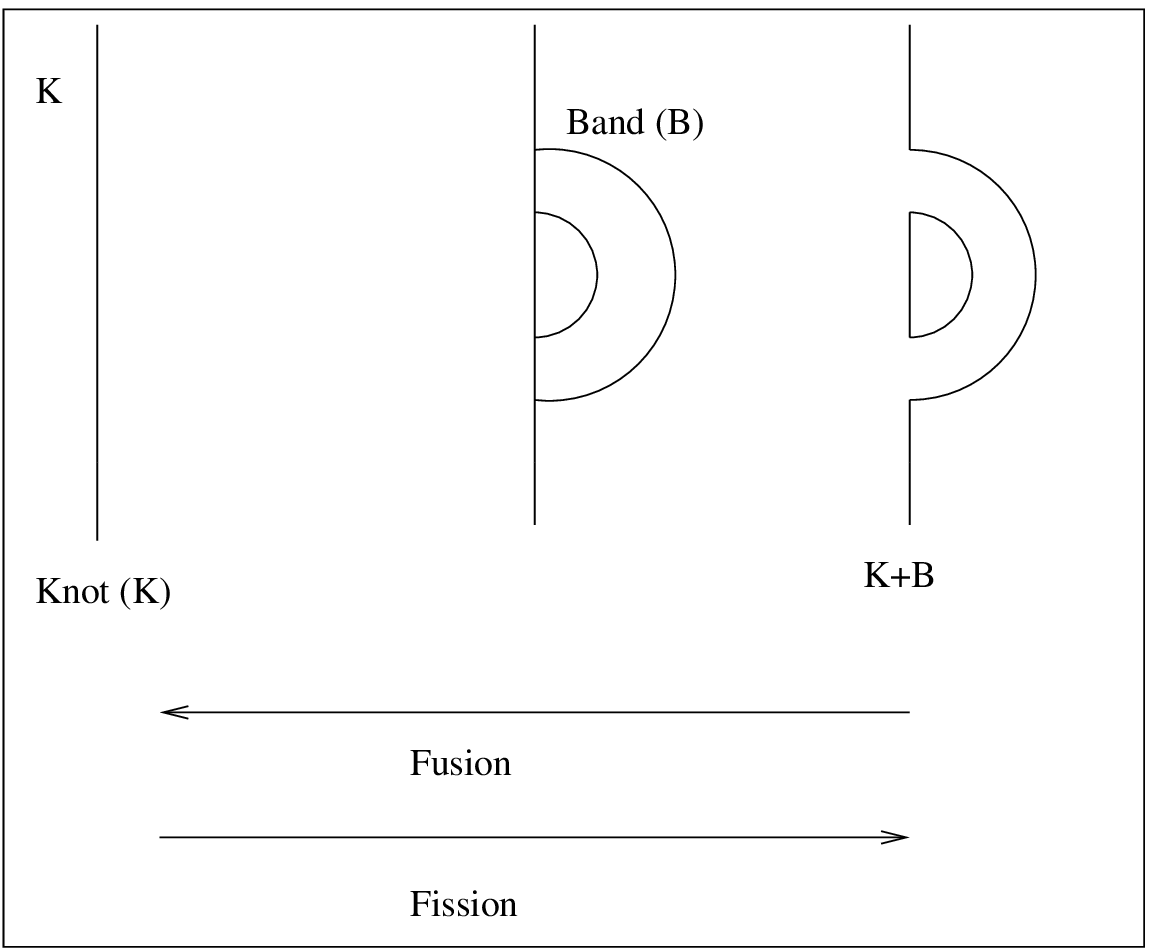}
\caption{Band}\label{fig14}
\end{center}
\end{figure}
\\
\begin{defn}\cite{Hoso}
Two links will be called link cobordic if one can be obtained from the other by a sequence of fusions and fissions.
This equivalence relation is denoted by $ \simeq$. $[L]$ denotes the link cobordism class of $L$.
\end{defn}
\begin{thm}\cite{Hoso}
The link cobordism classes with the conjunction operation form an abelian group, called the link cobordism
group and denoted by $ G_L$.
\end{thm}
\proof For two cobordism classes $[L_1]$ and $[L_2]$  the multiplication between them is well defined and given by
$$ [L_1] \& [L_2]=[L_1 \& L_2]$$
The zero of this operation is the class $[\bigcirc] $ which is the trivial link of a countable
number of components. The negative of $[L] $ is $ -[L]=[-L] $, where $-L$ denoted the reflected inverse
of $L$.
\endproof
\begin{lem}
For any link $L$, a conjunction $L\& -L$ is link cobordic to zero.
\end{lem}
To study the relation between the knot cobordism group $G_K$ and link cobordism group $G_L$ define a natural
mapping $ f:G_K \longrightarrow G_L $ which assigns to each knot cobordism class $[k]$ the corresponding link
cobordism class $[L]$ where $L$ is the knot $k$ viewed as a one-component link. We claim that $f$ is a homomorphism.
$f$ is well defined from the following lemma
\begin{lem}\cite{Hoso}\label{400}
Two knots are link cobordic if and only if they are knot cobordic.
\end{lem}
Now, $K_1 \# K_2$ is a fusion of $K_1\& K_2$ then $K_1 \# K_2$ is cobordic to $K_1\& K_2$, therefore
$f$ is a homomorphism. Again by using the lemma \ref{400}, if a knot is link cobordic to zero then it is also
knot cobordic to zero, and hence $ ker(f)$ consists of just $\bigcirc$.
\begin{lem}
f is an isomorphism of $G_K$ onto a subgroup of $G_L$.
\end{lem}
\begin{thm}\cite{Hoso}
$f(G_K)$ is a direct summand of $ G_L$ and it is a subgroup of $G_L$ whose elements have
total linking number zero. The other summand is isomorphic to the additive group of integers.
\end{thm}

\section{Graphs and cobordisms}
\subsection{Graph cobordism group}
In this section we construct cobordism groups for embedded graphs by extending the
notions of cobordisms used in the case of links.
\begin{defn}\label{GraphCobord1}
Two graphs $ E_1$ and $E_2$ are called cobordic if there is a surface $S$ have the
boundary $\partial S= E_1\cup -E_2$ with $E_1= S \cap (\s\times \{0\})$, $E_2=S \cap (\s\times \{1\})$
and we set $ E_1 \sim E_2$. Here by "surfaces" we mean
$2$-dimensional simplicial complexes that are PL-embedded in $\s \times [0,1]$.
$[E]$ denotes the cobordism class of the graph $E$.
\end{defn}
By using the graph cobordism classes and the conjunction operation $\&$, we can induce a graph cobordism group.
$ E_1 \& E_2 $ is a graph represented by the union of the two graphs $E_1 \cup E_2$ with mutually
disjoint 3-cells $ D_1$ and $ D_2$ containing (representatives of) $E_1$ and $ E_2$, respectively.
Here again we do not distinguish in the notation between the ambient isotopy classes
of embedded graphs (graph types) and a choice of representatives.
\begin{lem}\label{22332233}
The graph cobordism classes in the sense of Definition \ref{GraphCobord1}
with the conjunction operation form an abelian group called the graph cobordism
group and denoted by $G_E$.
\end{lem}
\proof
For two cobordism classes $[E_1]$ and $[E_2]$ the operation between them is given by
$$[E_1] \& [E_2]= [E_1 \& E_2].$$
This operation is well defined. To show that :
Suppose $ E_1 \sim F_1$, for two graphs $E_1$ and $F_1$. Then there exists a surface $S_1$
with boundary $ \partial S_1 = E_1\cup -F_1$. Suppose also, $ E_2 \sim F_2$, for two graphs $E_2$ and $F_2$.
Then there exists a surface $ S_2 $ with boundary $ \partial S_2= E_2\cup -F_2$. We want to show that
$ E_1 \& E_2 \sim F_1 \& F_2$, \ie we want to find a surface $S$ with boundary $\partial S= (E_1 \& E_2) \cup -(F_1 \& F_2)$.\\
Define the cobordism $S$ to be $ S_1 \& S_2$ where $S_1 \& S_2$ represents $ S_1 \cup S_2$ with mutually
disjoint 4-cells $D_1 \times [0,1]$ and $D_2 \times [0,1]$, containing $S_1$ and $S_2$ respectively with
$D_1 \times \{0\}$ containing $E_1$, $D_2 \times \{0\}$ containing $F_1$, $D_1 \times \{1\}$ containing $E_2$ and
$D_2 \times \{1\}$ containing $F_2$. The boundary of $S$ is given by,
$$\partial S= \partial (S_1\& S_2)=\partial S_1 \& \partial S_2=\partial S_1 \cup \partial S_2 =(E_1 \cup -F_1) \cup (E_2 \cup -F_2)=(E_1 \& E_2)\cup -(F_1 \& F_2) $$
Then the operation is well defined.
The zero of this operation is the class $[\includegraphics[width=0.5cm]{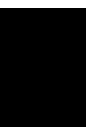}] $ which is the trivial graph of a countable
number of components. The negative of $[E] $ is $ -[E]=[-E] $, where $-E$ denoted the reflected inverse
of $E$.
\endproof
\subsection{Fusion and fission for embedded graphs}
We now describe a special kind of cobordisms between embedded graphs, namely the
basic cobordisms that correspond to attaching a 1-handle and that give rise to the analog
in the context of graphs of the operations of fusion and fission described already in the
case of links.
Let $E$ be a graph containing n-components with a split subgraph $E_s= G_1 \cup G_2 \cup G_3...\cup G_t$. We can define
a new graph $\hat{E}$ to be $ G_1 + G_2 + G_3...+G_t+\partial B_{t+1}+ \partial B_{t+2}+....+ \partial B_{n}$ where
$\{B_{t+1}, B_{t+2}, B_{t+3}....,B_{n}\}$ are disjoint bands in $\s$  spanning $E_s$. The operation $+$
means addition in the homology sense.
Put $ E_1 = E_s \cup  G_{t+1} \cup G_{t+2}....\cup G_{n}$ and
$E_2 = {\hat{E}} + G_{t+1} + G_{t+2}....+ G_{n}$. Now, The operation of replacing $E_1$ by $E_2$ is called {\em{\bf{fusion}}} and $E_2$ by $ E_1$ is called {\em{\bf{fission}}}.

\begin{figure}[htp]
\begin{center}
\includegraphics[width=10cm]{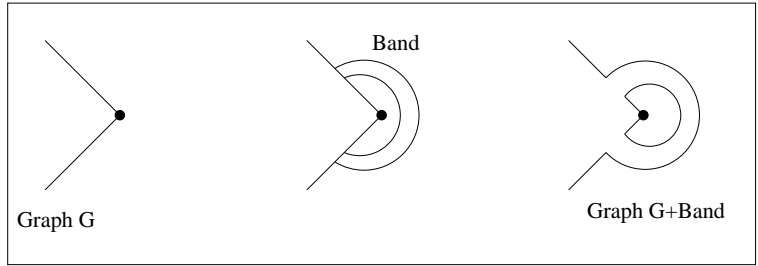}
\caption{Attaching Band to a Graph}\label{fig17.5}
\end{center}
\end{figure}

Notice that, in order to make sure that all resulting graphs will still have at least one vertex,
one needs to assume that the 1-handle is attached in such a way that there is at least an
intermediate vertex in between the two segments where the 1-handle is attached, as the
figure (\ref{fig17.5} )above.
\begin{rem}
Unlike the case of links, a fusion and fission for graphs does not necessarily change the number of components. For example see the figure (\ref{fig17}) below.
\end{rem}

\begin{figure}[htp]
\begin{center}
\includegraphics[width=10cm]{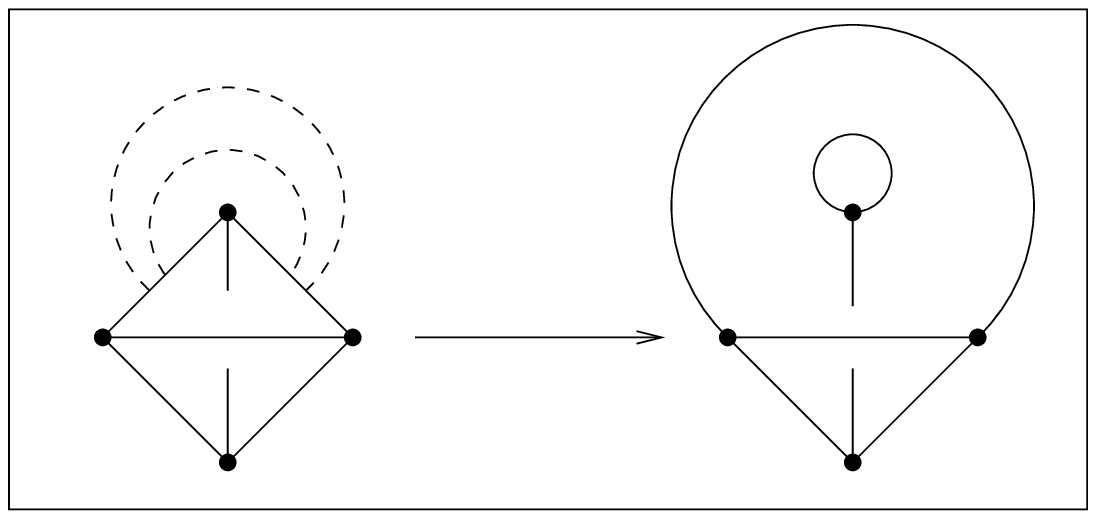}
\caption{Bands does not change the number of a graph components }\label{fig17}
\end{center}
\end{figure}

We can use the operations of fusion and fission described above to give another possible
definition of cobordism of embedded graphs.
\begin{defn}\label{GraphCobord2}
Two graphs will be called graph cobordic if one can be obtained from the other by a
sequence of fusions and fissions. We denote
this equivalence relation by $ \simeq$, and by $[E]$ the graph cobordism class of $E$.
\end{defn}
Thus we have two corresponding definitions for the graph cobordism group.
One can check from the definition of fusion and fission
that they gives the existence of a cobordism (surface) between two graphs $E_1$ and $E_2$.
\begin{lem}\label{Cobord21}
Two graphs $E_1$ and $E_2$ that are cobordant in the sense of Definition
\ref{GraphCobord2} are also cobordant in the sense of Definition \ref{GraphCobord1}.
The converse, however, is not necessary true.
\end{lem}
\proof As we have seen, a fusion/fission operation is equivalent to adding or removing
a band to a graph and this implies the existence of a saddle cobordism given by the attached
1-handle, as illustrated in figure (\ref{fig13}). By combining this cobordism with
the identity cobordism in the region outside where the 1-handle is attached, one obtains a
PL-cobordism between $E_1$ and $E_2$. This shows that cobordism in the sense of
Definition \ref{GraphCobord2}  implies cobordism in the sense of Definition  \ref{GraphCobord1}.
The reason why the converse need not be true is that, unlike what happens with
the cobordisms given by embedded smooth surfaces used in the case of links,
the cobordisms of graphs given by PL-embedded 2-complexes are not always
decomposable as a finite set of fundamental saddle cobordism given by a 1-handle.
Thus, having a PL-cobordism (surface in the sense of a 2-complex)
between two embedded graphs $E_1$ and $E_2$ does not necessarily
imply the existence of a finite sequence of fusions and fissions.
\endproof
\begin{lem}
The graph cobordism classes in the sense of Definition \ref{GraphCobord2} with the conjunction operation form an abelian group called the graph cobordism group and denoted by $G_F$.
\end{lem}
\proof
The proof is the same as the proof on lemma \ref{22332233} since fusion and
fission are a special case of cobordisms.
\endproof
The result of Lemma \ref{Cobord21} shows that there are different equivalence
classes $[E_1]\neq [E_2]$ in $G_F$ that are identified $[E_1]=[E_2]$ in $G_E$.
Thus, the number of cobordism classes when using Definition \ref{GraphCobord1}
is smaller that the number of classes by the fusion/fission method of
Definition \ref{GraphCobord2}.
\section{Khovanov Homology}

In the following we recall a homology theory for knots and links embedded in the
3-sphere. We discuss later  how to extend it to the case of embedded
graphs.
\\
\subsection{Khovanov Homology for links}
In recent years, many papers have appeared that discuss properties of Khovanov Homology theory,
which  was introduced in \cite{kh}. For each link  $ L \in \s $, Khovanov constructed
a bi-graded chain complex associated with the diagram $D$ for this link
$L$ and applied homology to get a group
$Kh^{i,j}(L)$, whose Euler characteristic is the normalized Jones polynomial.
$$ \sum_{i,j}(-1)^i q^j dim(Kh^{i,j}(L))=J(L)$$
He also proved that, given two diagrams $D$ and $D'$ for the same link, the corresponding
chain complexes are chain equivalent, hence, their homology groups are
isomorphic. Thus, Khovanov homology is a link invariant.

\subsection{The Link Cube}\label{licube}
\begin{figure}[htp]
\begin{center}
\includegraphics[width=6.5cm]{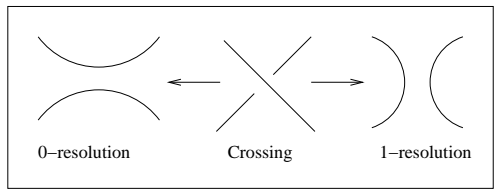}
\caption{0 and 1- resolutions to each crossing }\label{fig11}
\end{center}
\end{figure}

Let $L$ be a link with $n$ crossings. At any small neighborhood of a crossing we can
replace the crossing by a pair of parallel arcs and this operation is called a resolution.
There are two types of these resolutions called $0$-resolution (Horizontal resolution) and
$1$-resolution (Vertical resolution) as illustrated in figure (\ref{fig11}).\\
\\
We can construct a
$n$-dimensional cube by applying the $0$ and $1$-resolutions $n$ times to each crossing to
get $2^n$ pictures called smoothings (which are one dimensional manifolds) $S_\alpha$.
Each  of these can be indexed by a word $\mathbf{\alpha}$ of $n$ zeros and ones, \ie
$\mathbf{\alpha}  \in \{0,1\}^n$.
Let $\xi$ be an edge of the cube between two smoothings $S_{\alpha_{1}}$ and $S_{\alpha_{2}}$,
where $S_{\alpha_{1}}$ and $S_{\alpha_{2}}$ are identical smoothings except for a small
neighborhood around the crossing that changes from $0$ to $1$-resolution. To each edge
$\xi$ we can assign a cobordism $\Sigma_{\xi}$ (orientable surface whose boundary is the
union of the circles in the smoothing at either end)$$\Sigma_{\xi}: S_{\alpha_{1}} \longrightarrow S_{\alpha_{2}}$$
This $\Sigma_{\xi}$ is a product cobordism except in the neighborhood of the crossing, where it is  the obvious
saddle cobordism between the $0$ and $1$-resolutions.
Khovanov constructed a complex by applying a $1+1$-dimensional TQFT (Topological Quantum Field Theory)
which is a monoidal functor,
by replacing each vertex $S_\alpha$  by a graded vector space  $V_{\alpha}$
and each edge (cobordism) $\Sigma_{\xi}$ by a linear map $d_{\xi}$, and we set the group $CKh(D)$ to be the direct sum
of the graded vector spaces for all the vertices and the differential on the summand
$CKh(D)$ is a sum of the maps  $d_{\xi}$ for all edges $\xi$ such that Tail($\xi$)$=\alpha $ \ie
\begin{equation}\label{dv}
d^i(v)=\sum_{\xi} sign(-1)d_\xi(v)
\end{equation}
where $v \in V_\alpha \subseteq CKh(D)$ and $sign (-1)$  is chosen such that $d^2=0$.\\
An element of $CKh^{i,j}(D)$ is said to have homological grading $i$ and $q$-grading $j$ where
\begin{equation}\label{76}
i=|\alpha|-n_-
\end{equation}
\begin{equation}\label{77}
j=deg(v)+i+n_-+n_+
\end{equation}
for all $v \in V_\alpha \subseteq CKh^{i,j}(D)$, $|\alpha|$ is the number of 1's in $\alpha$, and
$n_-$, $n_+$ represent the number of negative and positive crossings respectively in the diagram $D$.

\subsection{Properties}\cite{pt1},\cite{lee}
Here we give some properties of Khovanov homology.
\begin{prop}\label{1234}
\begin{enumerate}
\item If $D'$ is a diagram obtained from $D$ by the application of a Reidemeister moves then
the complexes $ (CKh^{*,*}(D))$ and $ (CKh^{*,*}(D'))$ are homotopy equivalent.
\item For an oriented link $ L$ with diagram D, the graded Euler characteristic satisfies
\begin{equation}\label{12}
\sum (-1)^i qdim(Kh^{i,*}(L))=J(L)
\end{equation}
where $ J(L)$
is  the normalized Jones Polynomials for a link $L$ and $$\sum (-1)^i qdim(Kh^{i,*}(D))=\sum(-1)^i qdim(CKh^{i,*}(D))$$
\item Let $L_{odd}$ and $L_{even}$  be two links with odd and even number of components then $Kh^{*,even}(L_{odd})=0$ and
$Kh^{*,odd}(L_{even})=0$\\
\item
For two oriented link diagrams $D$ and $D'$, the chain complex of the disjoint union $ D \sqcup D'$ is given by
\begin{equation}\label{123}
CKh(D \sqcup D') =CKh(D) \otimes CKh (D').
\end{equation}
\item For two oriented links  $L$ and $L'$, the Khovanov homology of the disjoint union $ L \sqcup L'$ satisfies
$$Kh(L \sqcup L') = Kh(L) \otimes Kh (L').$$
\item Let $D$ be an oriented link diagram of a link $L$ with mirror image $ D^m$ diagram of the mirror link $ L^m$.
Then the chain complex $CKh(D^m)$ is isomorphic to the dual of $CKh(D)$ and $$ Kh(L) \cong Kh(L^m)$$
\end{enumerate}
\end{prop}

\subsection{Links Cobordisms}

\begin{figure}[htp]
\begin{center}
\includegraphics[width=5cm]{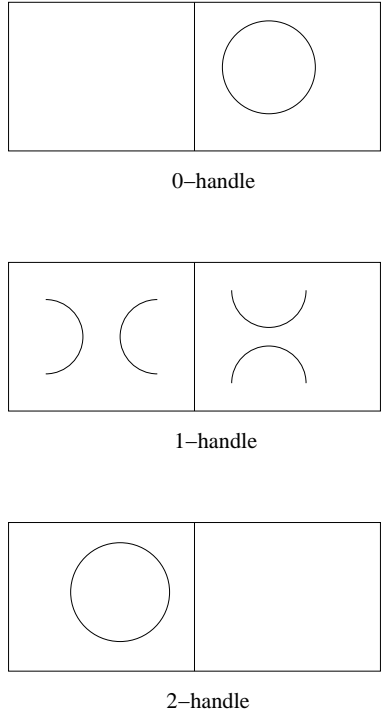}
\caption{Morse moves }\label{fig303}
\end{center}
\end{figure}

Let $ S $ $\subset \s \times [0,1]$ be a compact oriented surface
with boundary $\d S = L_0\cup -L_1$.
We assume that the $L_0$ and $L_1$ are links.
In this section we want to recall how one constructs a linear map between the
homologies of the boundary links by following Khovanov \cite{kh}.
The first idea is, we can decompose $S$ into elementary subcobordisms $ S_t$ for finitely many $ t \in [0,1]$ with
$$ S_t = S  \cap \s \times [0,t]$$ and $$ \d S_t = L_{t-1} \cup -L_{t}$$
where $L_{t-1}$ and $L_{t}$ are one dimensional manifolds, not necessary  links. Using a
small isotopy we can obtain that they are links for some $ t \in [0,1]$.
Here we assume that S is a smooth embedded surface.
A smooth embedded surface S can be represented by a one parameter family $D_t , t \in [0,1]$ of
planar diagrams
of oriented links $L_t$ for finitely many $ t \in [0,1]$  and this representation is called
a {\em movie} $M$. Between any two consecutive clips of a movie the
diagrams will differ by one of the ``Reidemeister moves" or ``Morse moves".
The Reidemeister moves are the first moves
in figure (\ref{fig21}) and the Morse moves are given in figure (\ref{fig303}).
\\
\\
These two types of moves will be called {\em local moves}. This means that between any two consecutive diagrams there is a local move either of
Reidemeister or of Morse type. The necessary condition is that the projection diagram  $D_0$ in the
first clip in $M$ should be the a projection of the link $ L_0$ and the projection diagram in the $ D_1$
in the last clip of the movie $M$ should be the projection of the link $L_1$ (boundary of $ S $).
Notice that the orientation
of $S$ induces an orientation on all intersection links $ L_t$. To show that, let $v$ be a tangent vector
to $L_t$. Then orient $v$ in the positive direction if $(v,w)$ gives the orientation  of $S$ where $w$ is the tangent vector to $S$ in the direction of increasing of $t$.
Khovanov constructed a chain map between complexes of two consecutive
diagrams that changed by a local move, hence a homomorphisms between
their homologies. The composition of these chain maps defines a homomorphism between
the homology groups of the diagrams of the boundary links.


\section{Homology theories for embedded graphs}
In this part we will present a method to extend Khovanov homology from links
to embedded graphs $G\subset \s$. Our construction is obtained by using Khovanov homology
for links, applied to a family of knots and links associated to an embedded graph. This
family is obtained by a result of  Kauffman \cite{kauff} as a fundamental topological
invariant of embedded graphs obtained by associating to an embedded graph $G$ in three-space a family of knots and links constructed by some operations of cutting graphs at vertices.
Before we give this construction, we motivate the problem of extending Khovanov homology
to embedded graphs by recalling another known construction of a homology theory, {\em graph homology}, which is defined for graphs. 

\subsection{Graph homology}
We recall here the construction and some basic properties of graph homology. As we discuss
below, graph homology can be regarded as a categorification of the chromatic polynomial
of a graph, in the same way as Khovanov homology gives a categorification of the
Jones polynomial of a link.
A construction is given for a graded homology theory for graphs whose graded Euler characteristic is the {\em Chromatic Polynomial}
of the graph \cite{laur}. Laure Helm-Guizon and Yongwu Rong used the same technique to get a graded chain complex. Their
construction depends on the edges in the vertices of the cube $\{0,1\}^n$ whose elements are connected subgraphs of the
graph $G$. In this subsection we recall the construction of Laure Helm-Guizon and Yongwu Rong.
\subsubsection{Chromatic Polynomial}
let $G$ be a graph with set of vertices $V(G)$ and set of edges $E(G)$. For a positive integer $t$, let
$\{ 1,2,...,t\}$ be the set of $t$-colors. A coloring of $G$ is an assignment of a $t$-color to each vertex
of $G$ such that vertices that are connected by an edge in $G$ always have different colors.
Let $P_G(t)$ be the number on $t$-coloring of $G$ \ie is the number of vertex colorings of $G$ with $t$
colors (in a vertex coloring two vertices are colored differently whenever they are connected by an edge $e$), then $P_G(t)$
satisfies the Deletion-Contraction relation
$$ P_G(t)=P_{G-e}(t)+P_{G/e}(t)$$
For an arbitrary edge $e \in E(G)$ we can define $ G-e$ to be the graph $G$ with deleted edge $e$, and
by $G/e$ the graph obtained by contracting edge $e$ \ie by identifying the vertices incident to $e$ and deleting $e$.
In addition to that $P_{K_n}(t)=t^n$ where $K_n$ is the graph with $n$ vertices and $n$ edges. $P_G(t)$ is called {\em Chromatic Polynomial}.
Another description can be give to $P_G(t)$, let $s \subset E(G)$, define $G_s$ to be the graph whose vertex set is the same
vertex set of $G$ with edge set $s$. Put $k(s)$ the number of connected components of  $G_s$. Then we have
$$ P_G(t)=\sum_{s\subseteq E(G)} (-1)^{|s|} t^{k(s)}$$

\subsubsection{Constructing $n$-cube for a Graph }
First we want to give an introduction to the type of algebra that we will use it in our work later.
\begin{defn}\cite{laur}
Let $\mathcal{V}=\oplus_i V_i$ be a graded $\mathbb{Z}$-module where $\{V_i\}$ denotes the set of homogenous
elements with degree $i$, and the graded dimension of $\mathcal{V}$ is the power series
$$ q dim\mathcal{V}=\sum_i q^i dim_{\mathbb{Q}}(V_i \otimes \mathbb{Q})$$.
\end{defn}
We can define the tensor product and directed sum for the graded $\mathbb{Z}$-module as follows:
\begin{thm}\cite{laur}
Let $\mathcal{V}$ and  $\mathcal{W}$ be a graded $\mathbb{Z}$-modules, then $\mathcal{V}\otimes \mathcal{W}$
and  $\mathcal{V}\oplus \mathcal{W}$ are both graded $\mathbb{Z}$-module with
\begin{enumerate}
\item $qdim(\mathcal{V}\oplus \mathcal{W})=qdim(\mathcal{V}) + qdim(\mathcal{W})$
\item $qdim(\mathcal{V}\otimes \mathcal{W})=qdim(\mathcal{V}) \cdot qdim(\mathcal{W})$
\end{enumerate}
\end{thm}

\begin{figure}[htp]
\begin{center}
\includegraphics[width=10cm]{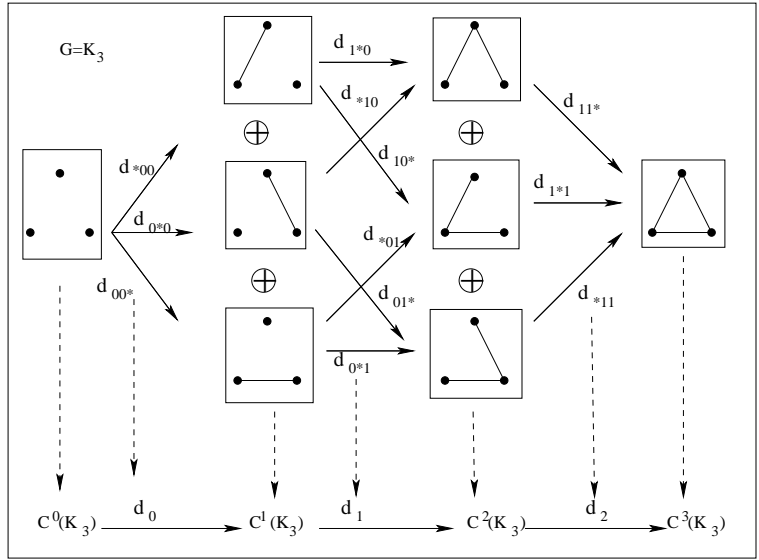}
\caption{$n$-dimensional cube for $K_3$ }\label{fig22}
\end{center}
\end{figure}

Let $G$ be a graph with edge set $E(G)$ and $n=|E(G)|$ represents the cardinality of $E(G)$.
We need first to order the edges in $E(G)$ and denote the edges by $\{e_1,e_2,...,e_n\}$. Consider the $n$-dimensional
cube $\{0,1\}^n$ \cite{laur},(see the figure (\ref{fig22})).
\\
\\
Each vertex can be indexed by a word $\alpha \in \{0,1\}^n$. This vertex $\alpha$
corresponded to a subset $s=s_{\alpha}$ of $ E(G)$. This is the set of edges of $G$
that are incident to the chosen vertex. Then
$e_i \in  s_{\alpha}$ if and only if  $\alpha_i=1$. Define $|\alpha|=\sum \alpha_i$
(height of $\alpha$) to be the number of $1's$ in $\alpha$ or equivalently the
number of edges in $s_{\alpha}$. We associate
to each vertex $\alpha$ in the cube $\{0,1\}^n$, a graded vector space $V_{\alpha}$
as follows \cite{laur}. Let $V_{\alpha}$ be a graded free $\mathbb{Z}$-module with
$1$ and $x$ basis elements with degree $0$ and $1$ respectively, then $V_{\alpha}=\mathbb{Z}\oplus \mathbb{Z}x$
with $qdim(V_{\alpha})=1+q$ and hence, $qdim(V_\alpha^{\otimes k})=(1+q)^k$.\\
Consider $G_{s_\alpha}$, the graph with vertex set $V(G)$ and edge set $s_\alpha$.
Replace each component of $G_{s_\alpha}$ by a copy of $V_{\alpha}$
and take the tensor product over all components.\\
Define the graded vector space $\cV_{\alpha}=V_{\alpha}^{\otimes k}$
where $k$ is the number of the components of $G_{s_\alpha}$.
Set the vector space $\cV$ to be the direct sum of the graded vector space for all the vertices.
The differential map $d^i$, defined by using the edges of the cube $ \{0,1\}^n$. We can label each
edge of $\{0,1\}^n$ by a sequence of $\{0,1,\star\}^n$ with exactly one $\star$. The tail of the edge
labeled by $\star=0$ and the head by $\star=1$. To define the differential we need first to define {\em Per-edge}
maps between the vertices of the cube $\{0,1\}^n$. These maps is defined to be a linear maps such that
every square in the cube $ \{0,1\}^n$ is commutative. Define the {\em per-edge} map $ d_{\xi}:\cV_{\alpha_1} \longrightarrow  \cV_{\alpha_2} $
for the edge $\xi$ with tail $\alpha_1$ and head $\alpha_2$ as follows: Take $\cV_{\alpha_i}=V^{\otimes k_i}$  for $i=1,2$
with $k_i$ is the number of the connected components of $G_{s_{\alpha_i}}$. Let $e$ be the edge and $s_{\alpha_2}=s_{\alpha_1} \cup \{e\}$, then there are two possible
cases.
First one (easy case): $d_{\xi}$ will be the identity map if the edge $e$ joins a component $r$ of $G_{s_{\alpha_1}}$
to itself. Then $k_1=k_2$ with a natural correspondence between the components
of $G_{s_{\alpha_1}}$ and $G_{s_{\alpha_2}}$.
Second one: if $e$ joins two  different components of $G_{s_{\alpha_1}}$, say $r_1$ and $r_2$, then $k_2=k_1-1$ and the components of $G_{s_{\alpha_2}}$
are $r_1\cup r_2 \cup \{e\}\cup....\cup r_{k_1}$. Define $d_{\xi}$ to be the identity map on the tensor
factor coming from $r_3,r_4,...,r_{k_1}$. Also define $ d_{\xi}$ on the remaining tensor factor to be the multiplication map
$V_{\alpha} \otimes V_{\alpha}  \longrightarrow V_{\alpha}$ sending $ x\otimes y $ to $xy$.
The differential $ d^i:\cV^i \longrightarrow \cV^{i+1}$ is given by
$$ d^i=\sum_{|\xi|=i} \sign(\xi)d_{\xi}$$
Where $\sign(\xi)$ is chosen so that $d^2=0$.
\begin{thm} \cite{pt1},\cite{laur}
The following properties hold for graph homology.
\begin{itemize}
\item The graded Euler characteristic for the graph homology given by
$$ \sum_{i,j} (-1)^i q^j dim(Kh^{i,j}(G))=P_G(t)$$
where $P_G(t)$ is the chromatic polynomial.

\item In graph homology a short exact sequence
$$ 0 \longrightarrow CKh^{i-1,j}(G/e) \longrightarrow CKh^{i,j}(G)\longrightarrow CKh^{i,j}(G-e)\longrightarrow 0 $$
can be constructed by using the deletion-contraction relation for a given edge $e \in G$.
This gives a long exact sequence
$$ \cdots \longrightarrow Kh^{i-1,j}(G/e)\longrightarrow Kh^{i,j}(G)\longrightarrow Kh^{i,j}(G-e)
\longrightarrow \cdots $$
\end{itemize}
\end{thm}

\subsection{Kauffman's invariant of Graphs}
\begin{figure}
\begin{center}
\includegraphics[width=8cm]{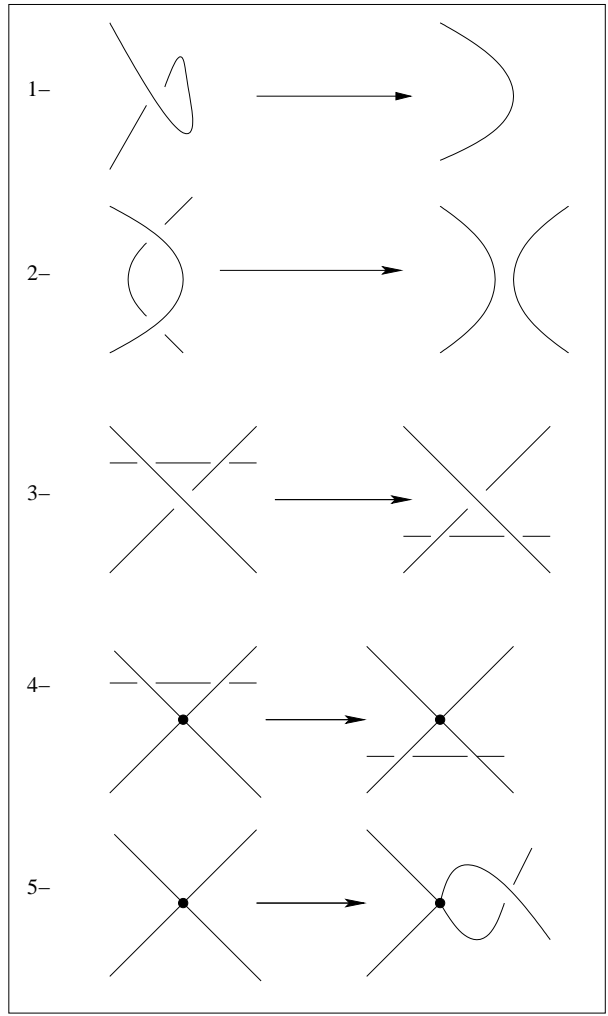}
\caption{Generalized Reidemeister moves by Kauffman}\label{fig21}
\end{center}
\end{figure}

\begin{figure}
\begin{center}
\includegraphics[width=8cm]{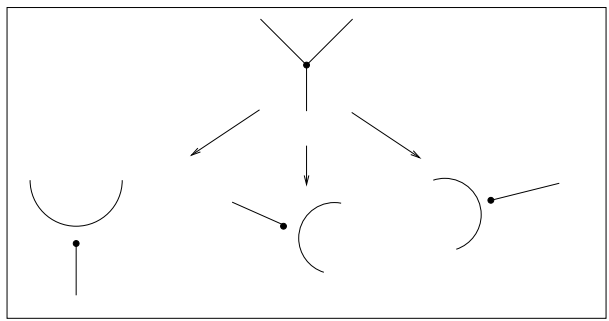}
\caption{local replacement to a vertex in the graph G}\label{fig18}
\end{center}
\end{figure}

\begin{figure}
\begin{center}
\includegraphics[width=8cm]{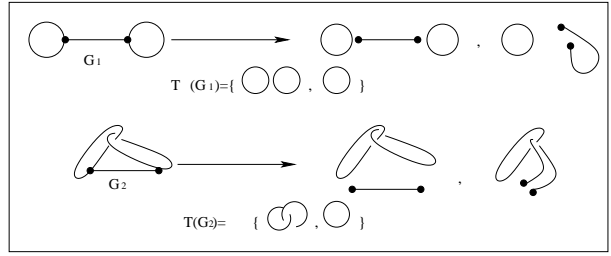}
\caption{Family of links associated to a graph}\label{fig19}
\end{center}
\end{figure}

We give now a survey of the Kauffman theory and show how to associate to an embedded
graph in $\s$ a family of knots and links. We then use these results to give our definition
of Khovanov homology for embedded graphs.
In \cite{kauff} Kauffman introduced a method for producing topological invariants of graphs
embedded in $\s$. The idea is to associate a collection of knots and links to a graph $G$ so that
this family is an invariant under the expanded Reidemeister moves defined by Kauffman and
reported here in figure (\ref{fig21}).
\\
\\
He defined in his work
an ambient isotopy for non-rigid (topological) vertices. (Physically, the rigid vertex
concept corresponds to
a network of rigid disks each with (four) flexible tubes or strings emanating from it.)
Kauffman proved that piecewise linear ambient isotopies of embedded graphs in $\s$
correspond to a sequence of generalized Reidemeister moves for planar diagrams of
the embedded graphs.
\begin{thm}\cite{kauff}
Piecewise linear (PL) ambient isotopy of embedded graphs is generated by the moves of figure (\ref{fig21}), that is, if two embedded graphs
are ambient isotopic, then any two diagrams of them are related by a finite sequence of the moves of figure (\ref{fig21}).
\end{thm}
Let $G$ be a graph embedded in $\s$. The procedure described by Kauffman of how to
associate to $G$ a family of  knots and links prescribes that we should make a local
replacement as in figure (\ref{fig18}) to each vertex in $G$.
Such a replacement at a vertex $v$ connects two edges and isolates all other edges at that vertex, leaving them as free ends. Let $r(G,v)$ denote the link formed by the closed curves formed by this process at a vertex $v$. One retains the link $r(G,v)$, while eliminating all the remaining unknotted arcs.
Define then $T(G)$ to be the family of the links $r(G,v)$ for all possible replacement choices,
$$ T(G)=\cup_{v\in V(G)} r(G,v). $$
For  example see figure (\ref{fig19}).
\\
\begin{thm}\cite{kauff}
Let $G$ be any graph embedded in $\s$, and presented diagrammatically.
Then the family of knots and links $T(G)$,
taken up to ambient isotopy, is a topological invariant of $G$.
\end{thm}

For example, in the figure (\ref{fig19}) the graph $G_2$ is not ambient isotopic to the graph $G_1$,
since $T(G_2)$ contains a non-trivial link.

\subsection{Definition of Khovanov homology for embedded graphs}
Now we are ready to speak about a new concept of Khovanov homology for embedded graphs by using
Khovanov homology for the links (knots) and Kauffman theory of associate a family of links
to an embedded graph $G$, as described above.
\begin{defn}
Let $G$ be an embedded graph with $ T(G)=\{ L_1,L_2,....,L_n\}$ the family of links associated to $G$
by the Kauffman procedure. Let $Kh(L_i)$ be the usual Khovanov homology of the link $L_i$ in this
family. Then the Khovanov homology
for the embedded graph $G$ is given by $$ Kh(G)=Kh(L_1) \oplus Kh(L_2)\oplus ....\oplus Kh(L_n)$$
Its graded Euler characteristic is the sum of the graded Euler characteristics of the
Khovanov homology of each link, \ie the sum of the Jones polynomials,
\begin{equation}
\sum_{i,j,k}(-1)^i q^j dim(Kh^{i,j}(L_k))=\sum_k J(L_k).
\end{equation}
\end{defn}
We show some simple explicit examples.
\begin{exa}
In figure (\ref{fig19}) $T(G_1)=\{\bigcirc \bigcirc ,\bigcirc\}$ then for $Kh(\bigcirc)=\mathbb{Q}$
$$ Kh(G_1)=Kh(\bigcirc \bigcirc) \oplus Kh(\bigcirc)$$
Now, from proposition \ref{1234} no.5
$$ Kh(G_1)= Kh(\bigcirc) \otimes Kh(\bigcirc) \oplus Kh(\bigcirc)$$
  $$ Kh(G_1)= \mathbb{Q} \otimes \mathbb{Q} \oplus \mathbb{Q}=\mathbb{Q} \oplus \mathbb{Q}$$
Another example comes from $T(G_2)=\{\includegraphics[width=0.6cm]{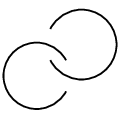}, \bigcirc\}$ then
$$Kh(G_2)=Kh(\includegraphics[width=0.6cm]{hopf.eps}) \oplus Kh(\bigcirc)$$
Since $ Kh^{0,0}(\bigcirc)=\mathbb{Q}$, and from \cite{pt1}

$$Kh(\includegraphics[width=0.6cm]{hopf.eps})= \begin{tabular}{|l|c|c|c|}\hline
\backslashbox{j}{i}  & -2 & -1 & 0 \\\hline
0 &    &      & $\mathbb{Q}$ \\\hline
-1 &    &      &    \\\hline
-2 &    &      & $\mathbb{Q}$    \\\hline
-3 &    &      &    \\\hline
-4 &$\mathbb{Q}$ &      &    \\\hline
-5 &    &      &    \\\hline
-6 &$\mathbb{Q}$ &      &    \\\hline
\end{tabular}
$$

Then,

$$Kh(G_2)= \begin{tabular}{|l|c|c|c|}\hline
\backslashbox{j}{i}  & -2 & -1 & 0 \\\hline
0 &    &      & $\mathbb{Q}\oplus  \mathbb{Q}$ \\\hline
-1 &    &      &    \\\hline
-2 &    &      & $\mathbb{Q}$    \\\hline
-3 &    &      &    \\\hline
-4 &$\mathbb{Q}$ &      &    \\\hline
-5 &    &      &    \\\hline
-6 &$\mathbb{Q}$ &      &    \\\hline
\end{tabular}
$$
\end{exa}

\bibliographystyle{amsplain}

\providecommand{\bysame}{\leavevmode\hbox to3em{\hrulefill}\thinspace}
\providecommand{\MR}{\relax\ifhmode\unskip\space\fi MR }
\providecommand{\MRhref}[2]{%
  \href{http://www.ams.org/mathscinet-getitem?mr=#1}{#2}
}
\providecommand{\href}[2]{#2}
\begin{thebibliography}{}

\end{thebibliography}


\begin{thebibliography}{33333}
\bibitem{Ba}  D.~Bar-Natan, {\em Khovanov's homology for tangles
and cobordisms}, Geometry and Topology, Vol.9 (2005) 1443--1499.
\bibitem{foxmi} R.H.Fox and J.W. Milnor {\em singularity of 2-spheres in 4-space and equivalence knots}
(Abstract), Bull.Amer. Math. Soc.,Vol 63 (1957), pp 406
\bibitem{fox} R.H.~Fox, {\em Covering spaces with singularities},
in ``Algebraic Geometry and Topology. A symposium in honour of
S.Lefschetz'' Princeton University Press 1957, 243--257.
\bibitem{fox2} R.H.~Fox, {\em A quick trip through knot theory}. in
``Topology of 3-manifolds and related topics'' Prentice-Hall, 1962,
pp. 120--167.
\bibitem{laur} Laure Helm-Guizon and Yongwu Rong  {\em Graph Cohomologies from Arbitrary Algebras},
mathQA/0506023v1, (2005).
\bibitem{Hoso} F.~Hosokawa, {\em A concept of cobordism between links},
Ann. Math. Vol.86 (1967) N.2, 362--373.
\bibitem{Jaco} M. Jacobsson {\em An Invariant of Link Cobordisms from Khovanov Homology}, Algebraic and Geometry Topology,
Vol.4, (2004), 1211-1251
\bibitem{kauff} Louis H.Kauffman {\em Invariants of Graphs in Three-Space}, Transactions of the American Mathematical
Society, Vol.311, No.2, (1989).
\bibitem{kh} M.~Khovanov, {\em A categorification of the Jones
polynomial}. Duke Math. J.  101  (2000),  no. 3, 359--426.
\bibitem{lee} E. S.Lee  {\em An Endomorphism of the Khovanov Invariant}, math.GT/0210213v3, (2004)

\bibitem{jr} Jacob Rasmussen {\em Khovanov Homology and the Slice Genus}, math.GT/0402131v1,(2004)
\bibitem{Tan} K.Taniyama, {\em Cobordism, homotopy and homology of graphs in $\R^3$},
Topology, Vol.33 (1994) 509-523.

\bibitem{TB} C.H.~Taubes, J.~Bryan, {\em Donaldson-Floer theory}, in ``Gauge
theory and the topology of four-manifolds'' (Park City, UT, 1994),  195--221,
IAS/Park City Math. Ser., 4, Amer. Math. Soc., Providence, RI, 1998.
\bibitem{pt1} Paul Turner {\em Five Lectures on Khovanov Homology}, math.GT/0606464,(2006)

\end{thebibliography}

\end{document}